\documentclass[draft]{agujournal2019}
\usepackage{url} 
\usepackage{lineno}
\usepackage[finalnew]{trackchanges} 
\usepackage{soul}

\usepackage{bm}
\usepackage{upgreek}
\usepackage{amsmath}
\usepackage{amssymb} 
\usepackage{amsfonts}
\usepackage{algorithm}
\usepackage{algorithmic}
\usepackage{color}
\usepackage{ragged2e}

\justifying

\draftfalse

\journalname{Water Resources Research}

\begin{document}

%
%


\title{Using Deep Learning to Improve Ensemble Smoother: Applications to Subsurface Characterization}

%
%




\authors{Jiangjiang Zhang\affil{1}, Qiang Zheng\affil{2}, Laosheng Wu\affil{3}, and Lingzao Zeng\affil{2}}


\affiliation{1}{Yangtze Institute for Conservation and Development, Hohai University, Nanjing, China,}
\affiliation{2}{Zhejiang Provincial Key Laboratory of Agricultural Resources and Environment, Institute of Soil and Water Resources and Environmental Science, College of Environmental and Resource Sciences, Zhejiang University, Hangzhou, China,}
\affiliation{3}{Department of Environmental Sciences, University of California, Riverside, California, USA.}




\correspondingauthor{L. Zeng}{lingzao@zju.edu.cn}




\begin{keypoints}
\item Ensemble smoother using the Kalman formula is constrained by the Gaussian assumption
\item We use deep learning to formulate a new variant of ensemble smoother
\item The new method can produce better results in problems involving non-Gaussian distributions
\end{keypoints}

%
%

%
%


\begin{abstract}
Ensemble smoother (ES) has been widely used in various research fields to reduce the uncertainty of the system-of-interest. However, the commonly-adopted ES method that employs the Kalman formula, that is, ES$_\text{(K)}$, does not perform well when the probability distributions involved are non-Gaussian. To address this issue, we suggest to use deep learning (DL) to derive an alternative analysis scheme for ES in non-Gaussian data assimilation problems. Here we show that the DL-based ES method, that is, ES$_\text{(DL)}$, is more general and flexible. In this new scheme, a high volume of training data is generated from a relatively small-sized ensemble of model parameters and simulation outputs, and possible non-Gaussian features can be preserved in the training data and captured by an adequate DL model. This new variant of ES is tested in two subsurface characterization problems with or without the Gaussian assumption. Results indicate that ES$_\text{(DL)}$ can produce similar (in the Gaussian case) or even better (in the non-Gaussian case) results compared to those from ES$_\text{(K)}$. The success of ES$_\text{(DL)}$ comes from the power of DL in extracting complex (including non-Gaussian) features and learning nonlinear relationships from massive amounts of training data. Although in this work we only apply the ES$_\text{(DL)}$ method in parameter estimation problems, the proposed idea can be conveniently extended to analysis of model structural uncertainty and state estimation in real-time forecasting problems. \end{abstract}


%
%

%


%
%
%
%

\section{Introduction}

Numerical models have been widely used in science and engineering to gain a better understanding of the concerned process(es), and to help hypothesis testing and decision making. In many research fields of geosciences, complexity of the system-of-interest makes accurate predictions in both the space and time domain very challenging \cite{baartman2020,kavetski2006a,kavetski2006b,refsgaard2012,ruddell2019}. This is largely caused by our incomplete knowledge and insufficient observations of the system. To improve our predictive ability and scientific understanding of the system, it is important to combine the numerical model (i.e., the theory) with observations (i.e., the data), which can be realized through data assimilation \cite<DA;>[]{carrassi2018,evensen2009}. DA is usually carried out in the following way. At any update time, one first makes a forecast from the background information. The forecast variables can be initial/boundary conditions, parameters, simulation outputs, model errors, or their combinations \cite{carrassi2018,chen2006,dechant2011,evensen2009,evensen2019,wang2020robust,xue2014,zhang2019}. Then one calculates the difference (which is usually called the innovation) between the observations and the corresponding model outputs mapped from the forecast with a linear/nonlinear operator. The innovation vector provides new information, based on which some update (or correction) to the forecast can be made.

Theoretically, one can view DA from a Bayesian perspective. However, fully Bayesian DA methods, such as particle filter \cite<PF;>[]{doucet2000,moradkhani2005}, although theoretically complete, can be computationally prohibitive for high-dimensional problems. When the probability distributions are assumed to be Gaussian, that is, the update is only based on the mean and (co)variance, DA can be implemented efficiently. In linear systems, the Kalman filter is the optimal DA method \cite{kalman1960}. For nonlinear dynamics, one can linearize the state equations and apply the extended Kalman filter \cite<EKF;>[]{gelb1974}. In EKF, the analysis scheme is based on the Jacobian of the state equations. When the problem at hand is high-dimensional and nonlinear, the performance of EKF becomes unsatisfactory, then the ensemble Kalman filter (EnKF) proposed by \citeA{evensen1994} can be adopted as a promising alternative. EnKF is a Monte Carlo implementation of the Kalman filter to perform sequential DA. When the purpose is parameter estimation, it will be more convenient to perform a global update using the entire historical observations, instead of applying the sequential analysis scheme of EnKF. In this case, ensemble smoother \cite<ES;>[]{vanleeuwen1996} can be employed. It has been shown that, ES can obtain similar results to EnKF, but with a much lower computational cost \cite{li2018,skjervheim2011}. When the system is highly nonlinear, iterative applications of both EnKF and ES are needed \cite{chen2012,emerick2012,emerick2013,gu2007,lorentzen2011}. In the past decades, EnKF and its variants have been extensively used in various research fields, for example, meteorology \cite{houtekamer2016}, oceanography \cite{chen2009,simon2009}, hydrology \cite{chen2006,reichle2008,schoniger2012,xie2010}, and petroleum engineering \cite{aanonsen2009,emerick2012,gu2007}, just to name a few .

Nevertheless, in many situations, distributions of parameters, simulation outputs, measurement errors and so forth can be obviously non-Gaussian \cite{mandel2009,schoups2010,sun2009,zhou2011}. Then the direct use of a Kalman-based DA method, such as EnKF or ES, becomes inappropriate. To address this issue, different strategies have been proposed. For example, when the probability distribution is multi-modal, one can first turn the forecast ensemble into several clusters, and each cluster can be roughly approximated by a Gaussian distribution and updated with a Kalman-based DA method \cite{bengtsson2003,dovera2011,elsheikh2013,sun2009,zhang2018}; Some other researchers suggested to re-parameterize non-Gaussian variables (e.g., conductivity field in a channelized aquifer) to be Gaussian distributed with anamorphosis function \cite{schoniger2012,simon2009}, level set \cite{chang2010}, or normal-score transformation \cite{li2011,li2018,xu2016,zhou2011}, etc., so that EnKF or its variants can be properly implemented; Another kind of approaches first use a Kalman-based DA method to update the forecast ensemble, and then use, for example, multiple-point geostatistics \cite{cao2018,jafarpour2011,kang2019,sarma2008}, or a more general DA method like PF \cite{mandel2009}, to reconstruct the non-Gaussian target distributions. Note that although the above-mentioned strategies worked well in different applications, they have not changed the DA methods that were used. In other words, the DA methods are still more-or-less constrained by the Gaussian assumption, and these strategies either apply some pre-treatment to fulfill this assumption, or use some post-treatment to fix the DA results.

In this work, we propose a new DA method that is free from the Gaussian assumption, and at the same time is computationally feasible for high-dimensional problems. Before introducing the basic idea behind this method, let's first go back to the general process of DA that has been demonstrated earlier. Essentially, DA works by updating (or correcting) the forecast from the innovation (i.e., the difference between observations and the corresponding model outputs). In the Kalman-based DA methods, a linear mapping from the innovation vector to the update vector is calculated from the forecast states based on the Kalman formula. This mapping, usually called the Kalman gain, only uses the first and second-order statistical moments. Then it is natural to wonder whether we can obtain a more general (i.e., free from the Gaussian assumption), and possibly nonlinear mapping to update the forecast states. In the past years, machine learning, especially deep learning (DL), has been extensively used in different fields, including hydrology and water resources, to extract complex features and learn nonlinear relationships from data \cite{goodfellow2016,lecun2015,shen2018,shen2018hess}. It has motivated us to reformulate DA, especially the Kalman-based DA methods, through obtaining a possibly nonlinear mapping from the innovation vector to the update vector with DL. Now one question arises, that is, how can we generate a high volume of training data that is usually required to feed a DL model, when a large number of system model evaluations are not affordable? To address this issue, we come up with a simple solution. In the forecast ensemble with $N_\text{e}$ samples of model parameters and simulation outputs, if we pick out one arbitrary sample as the hypothetical truth, and generate synthetic observations by perturbing the ``true" model outputs with random errors, we can obtain $N_\text{e}-1$ pairs of innovation and update vectors. It means that we pick two elements out of the $N_\text{e}$ samples at a time without repetition, where one element is regarded as the hypothetical truth. From the basic theory of combination, we can generate training data with $C(N_\text{e},2)=N_\text{e}(N_\text{e}-1)/2$ unique samples for DL. In the $C(N_\text{e},2)$ samples, non-Gaussian features of model parameters and measurement data can be preserved in the synthetic innovation and update vectors, and captured by an adequate DL model. Finally, we can input the actual $N_\text{e}$ innovation vectors (calculated from the measurements and the forecast ensemble) to the DL-based mapping, and obtain $N_\text{e}$ update vectors to correct our forecast. There are several advantages in using DL to derive the new mapping, that is, 1) this mapping can be nonlinear, and complex (e.g., non-Gaussian) features can be extracted; and 2) one can avoid the calculation of Jacobian (in EKF) and/or covariance matrices, as well as the inversion of matrix.

In this work, to verify the validity of the proposed idea, we use DL to improve ES, and test the resulting ES$_\text{(DL)}$ method in two subsurface characterization problems with or without the Gaussian assumption. In subsurface characterization, DL has been used as a powerful tool to address a wide range of challenges. For example, DL can be used to effectively reduce the dimensionality of model parameters \cite{laloy2017}, or quickly generate random realizations of (e.g., non-Gaussian) geological formations from training data \cite<e.g., image data;>[]{laloy2018}, both of which can improve the performance of geostatistical inversion; To alleviate the high computational cost caused by repetitive evaluations of complex, high-dimensional groundwater models, DL-based surrogates have been built and used in uncertainty quantification and DA for subsurface systems \cite{mo2019,mo2020,tripathy2018}; In \citeA{mo2020}, a Kalman-based DA method proposed in our earlier work \cite<ILUES;>[]{zhang2018} was adopted, and another DL model was used to parameterize non-Gaussian conductivity field; If physical laws are preserved when training a DL model, a better performance of the resulting surrogate can be obtained \cite{wang2020}; Besides constructing a forward mapping from parameters to simulation outputs as a surrogate model, \citeA{sun2018} also identified a reverse mapping from simulation outputs to parameters. For more applications of DL in hydrology and water resources, one can refer to \cite{shen2018,shen2018hess}.

The rest of this paper is organized as follows. In section 2, we first introduce how to implement ES using the Kalman formula, that is, ES$_\text{(K)}$, to estimate unknown model parameters from indirect measurement data. In light of the limitations of ES$_\text{(K)}$, we then propose a more general method, ES$_\text{(DL)}$, that uses DL to extract non-Gaussian features and learn a nonlinear analysis scheme. To verify the performance of ES$_\text{(DL)}$, two cases of subsurface characterization with or without the Gaussian assumption are tested in section 3. Here we are concerned with benchmarking analysis of
the two ES methods. Finally, in section 4, we conclude this paper and discuss the pros and cons of the new method.

\section{Methods}

Let's assume that the system-of-interest is simulated by a numerical model, $f(\cdot)$, and this process can be expressed in the following compact form,
\begin{linenomath*}
	\begin{equation}
	\widetilde{\textbf{y}}=f(\textbf{m})+\upepsilon,
	\label{eq: numerical_model}
	\end{equation}
\end{linenomath*}
where $\widetilde{\textbf{y}}\in \mathbb{R}^{N_\text{y}}$ are observations of the system, $\textbf{m}\in \mathbb{R}^{N_\text{m}}$ are the model parameters, and $\upepsilon\in\mathbb{R}^{N_\text{y}}$ are the error term. In subsurface characterization, the model parameters include the spatial/temporal distributions of contaminants and/or subsurface properties, which are generally difficult or even impossible to be measured directly. Meanwhile, observations of some state variables, such as hydraulic head, solute concentration, temperature, and electromagnetic signals, can be monitored continuously and affordably. These observations, that is, $\widetilde{\textbf{y}}$, contain information about the unknown model parameters, $\textbf{m}$. To improve our knowledge of $\textbf{m}$, we can perform data assimilation conditioned on these measurement data.

\subsection{Ensemble Smoother Using the Kalman Formula: ES$_\text{(K)}$}

As an efficient and robust data assimilation method, EnKF has been extensively used in various research fields to reduce the uncertainty of the system-of-interest \cite{evensen2009}. When one's purpose is parameter estimation, a variant of EnKF, that is, ES, can be adopted as a suitable method \cite{vanleeuwen1996}. Below we will introduce how to implement ES that uses the Kalman formula, that is, ES$_\text{(K)}$, to estimate unknown model parameters, $\textbf{m}$, from indirect measurement data, $\widetilde{\textbf{y}}$.

Here, a prior distribution, $p(\textbf{m})$, is used to represent our background knowledge of the values of $\textbf{m}$. From $p(\textbf{m})$ we can draw $N_\text{e}$ random samples to form the forecast (or prior) ensemble, that is, $\textbf{M}^{(0)}=\{\textbf{m}_1^{(0)},...,\textbf{m}_{N_\text{e}}^{(0)}\}$. Then we calculate the corresponding model outputs by running the numerical model, that is, $\textbf{Y}^{(0)}=\{f(\textbf{m}_1^{(0)}),...,f(\textbf{m}_{N_\text{e}}^{(0)})\}$. Using the ES$_\text{(K)}$ method, we can update each sample, $\textbf{m}_i^{(0)}$, $i=1,...,N_\text{e}$, in the forecast ensemble, $\textbf{M}^{(0)}$, conditioned on the measurement data, $\widetilde{\textbf{y}}$,
\begin{linenomath*}
	\begin{equation}
	\textbf{m}_i^{(1)}=\textbf{m}_i^{(0)}+\textbf{C}_{\text{MY}}^{(0)}\left(\textbf{C}_{\text{YY}}^{(0)}+\textbf{R}\right)^{-1}\left[\widetilde{\textbf{y}}+\upepsilon_i-f(\textbf{m}_i^{(0)})\right],
	\label{eq: ES_K}
	\end{equation}
\end{linenomath*} 
where $\textbf{M}^{(1)}=\{\textbf{m}_1^{(1)},...,\textbf{m}_{N_\text{e}}^{(1)}\}$ is the updated ensemble, $\textbf{C}_{\text{MY}}^{(0)}$ is the cross-covariance matrix between model parameters and simulation outputs (calculated from $\textbf{M}^{(0)}$ and $\textbf{Y}^{(0)}$), $\textbf{C}_{\text{YY}}^{(0)}$ is the auto-covariance matrix of model outputs (calculated from $\textbf{Y}^{(0)}$), and $\upepsilon_i$ is a random realization of measurement error with covariance $\textbf{R}$. In ES$_\text{(K)}$, the analysis scheme is essentially linear, and the distributions of model parameters and measurement data should be close to Gaussian. At this point the applicability of ES$_\text{(K)}$ is limited.

\subsection{Using Deep Learning to Improve Ensemble Smoother: ES$_\text{(DL)}$}

In the past few years, DL has been extensively used to learn complex patterns and nonlinear relationships from data. The general applicability of DL has motivated us to reformulate the analysis scheme of ES to make it more capable. The new ES method is termed ES$_\text{(DL)}$. Nowadays, a plenty of powerful DL architectures have been proposed by the machine learning communities. Here we are only left with choosing which relationship to learn and how to generate enough data to train the DL model. As the theory of DL itself is not the focus of this work, we decide not to provide the details here. Interested readers are suggested to refer to \cite{goodfellow2016,lecun2015}. Moreover, architectures of the DL models used in this work will be given in section 3 (Figures \ref{fig: ex1_dl} and \ref{fig: ex2_dl}).

Let's rewrite equation (\ref{eq: ES_K}) in a more general form,
\begin{linenomath*}
	\begin{equation}
	\Delta\textbf{m}_i=\mathcal{G}(\Delta\textbf{y}_i),
	\label{eq: ES_general}
	\end{equation} 
\end{linenomath*}
where $\Delta\textbf{m}_i=\textbf{m}_i^{(1)}-\textbf{m}_i^{(0)}$ is the update vector, $\Delta\textbf{y}_i=\widetilde{\textbf{y}}+\upepsilon_i-f(\textbf{m}_i^{(0)})$ is the innovation vector, and $\mathcal{G}$ is a mapping from $\Delta\textbf{y}_i$ to $\Delta\textbf{m}_i$. In ES$_{(\text{K})}$, $\mathcal{G}$ is defined by the Kalman gain matrix, $\textbf{K}=\textbf{C}_{\text{MY}}^{(0)}\left(\textbf{C}_{\text{YY}}^{(0)}+\textbf{R}\right)^{-1}$. Thus, the relationship between $\Delta\textbf{y}_i$ and $\Delta\textbf{m}_i$ is linear. Here we aim to use DL to derive a possibly nonlinear mapping, $\mathcal{G}_{\text{DL}}$, from $\Delta\textbf{y}_i$ to $\Delta\textbf{m}_i$.

It is evident that the input data of $\mathcal{G}_{\text{DL}}$ are the displacement vector in model simulations (corrupted by some errors), and the outputs are the corresponding distance in the parameter space. Based on this finding, we can generate a high volume of training data from the forecast ensemble, $\{\textbf{M}^{(0)}~\textbf{Y}^{(0)}\}$. In $\{\textbf{M}^{(0)}~\textbf{Y}^{(0)}\}$, if we take two samples at a time without repetition, there will be $N=N_\text{e}(N_\text{e}-1)/2$ combinations in total. That is to say, we can obtain $\textbf{D}_\text{in}^{(0)}=\{f(\textbf{m}_i^{(0)})-f(\textbf{m}_j^{(0)})+\upepsilon_{ij}|i=1,...,N_\text{e}-1,~i<j\le N_\text{e}
\}$ as the inputs to the DL model, and $\textbf{D}_\text{out}^{(0)}=\{\textbf{m}_i^{(0)}-\textbf{m}_j^{(0)}|i=1,...,N_\text{e}-1,~i<j\le N_\text{e}\}$ as the output data, where $\upepsilon_{ij}$ are random realizations of the measurement error. Here, non-Gaussian features in the model parameters and observations can be preserved in the training data, $\textbf{D}^{(0)}=\{\textbf{D}_\text{in}^{(0)}~\textbf{D}_\text{out}^{(0)}\}$, and captured by the DL model. When the evaluation of $f(\cdot)$ is time-consuming, one usually can only afford a limited number of model runs (e.g., $N_\text{e}=200$). In this case, the number of samples in $\textbf{D}^{(0)}$ is still considerable ($N=19,900$).

After training, a possibly nonlinear mapping $\mathcal{G}_\text{DL}$ can be obtained. Then we can use $\mathcal{G}_\text{DL}$ to update each sample, $\textbf{m}_i^{(0)}$, $i=1,...,N_\text{e}$, in the forecast ensemble, $\textbf{M}^{(0)}$, conditioned on the measurement data, $\widetilde{\textbf{y}}$,
\begin{linenomath*}
	\begin{equation}
	\textbf{m}_i^{(1)}=\textbf{m}_i^{(0)}+\mathcal{G}_\text{DL}\left(\widetilde{\textbf{y}}+\upepsilon_i-f(\textbf{m}_i^{(0)})\right).
	\label{eq: ES_DL}
	\end{equation}
\end{linenomath*}
Then we use the updated ensemble, $\textbf{M}^{(1)}=\{\textbf{m}_1^{(1)},...,\textbf{m}_{N_\text{e}}^{(1)}\}$, to represent our new knowledge about the model parameters. \add[editor]{\protect It is noted here that although ES$_\text{(DL)}$ is developed from ES$_\text{(K)}$, the two methods differ from each other theoretically. ES$_\text{(K)}$ can be derived from the least-squares approach \cite{anderson2003}. However, in ES$_\text{(DL)}$, the analysis scheme is directly learned from training data of (synthetic) innovation and update vectors with an adequate DL model.}

For highly nonlinear problems, one single update of the model parameters with ES$_\text{(K)}$ or ES$_\text{(DL)}$ may not be sufficient. Here, we suggest to adopt the multiple data assimilation scheme proposed by \citeA{emerick2013} to address this issue. In this scheme, the measurement data are assimilated $N_\text{iter}$ times. To make sure that the finally obtained results are reasonable, the measurement error (including the corresponding covariance matrix $\textbf{R}$, if used) should be inflated by a factor of $\alpha_t$ (for $\textbf{R}$ the factor is the square of $\alpha_t$) in iteration $t$, $t=1,...,N_\text{iter}$. The factors should satisfy $\sum_{t=1}^{N_\text{iter}}1/(\alpha_t)^2=1$, and a convenient choice is $\alpha_t=\sqrt{N_\text{iter}}$. The analysis scheme for ES$_{\text{(K)}}$ becomes,
\begin{linenomath*}
	\begin{equation}
	\textbf{m}_i^{(t)}=\textbf{m}_i^{(t-1)}+\textbf{C}_{\text{MY}}^{(t-1)}\left[\textbf{C}_{\text{YY}}^{(t-1)}+(\alpha_{t})^2\textbf{R}\right]^{-1}\left[\widetilde{\textbf{y}}+\alpha_t\upepsilon_i-f(\textbf{m}_i^{(t-1)})\right].
	\label{eq: ES_K_MDA}
	\end{equation}
\end{linenomath*}
In ES$_\text{(DL)}$, we first generate training data from $\{\textbf{M}^{(t-1)}~\textbf{Y}^{(t-1)}\}$ as, $\textbf{D}_\text{in}^{(t-1)}=\{f(\textbf{m}_i^{(t-1)})-f(\textbf{m}_j^{(t-1)})+\alpha_t\upepsilon_{ij}|i=1,...,N_\text{e}-1,~i<j\le N_\text{e}\}$, and $\textbf{D}_\text{out}^{(t-1)}=\{\textbf{m}_i^{(t-1)}-\textbf{m}_j^{(t-1)}|i=1,...,N_\text{e}-1,~i<j\le N_\text{e}\}$. Then a mapping, $\mathcal{G}_\text{DL}^{(t-1)}$, is discovered from $\{\textbf{D}_\text{in}^{(t-1)}~\textbf{D}_\text{out}^{(t-1)}\}$ with DL. Finally, each sample in the forecast ensemble is updated as,
\begin{linenomath*}
	\begin{equation}
	\textbf{m}_i^{(t)}=\textbf{m}_i^{(t-1)}+\mathcal{G}_\text{DL}^{(t-1)}\left(\widetilde{\textbf{y}}+\alpha_{t}\upepsilon_i-f(\textbf{m}_i^{(t-1)})\right).
	\label{eq: ES_DL_MDA}
	\end{equation}
\end{linenomath*}
For both ES$_\text{(K)}$ and ES$_\text{(DL)}$, our final knowledge of the model parameters is represented by $\textbf{M}^{(N_\text{iter})}=\{\textbf{m}_1^{(N_\text{iter})},...,\textbf{m}_{N_\text{e}}^{(N_\text{iter})}\}$.

\section{Illustrative Case Studies}
\subsection{Example 1: A Gaussian Case}

\begin{figure}
	\noindent\includegraphics[width=\textwidth]{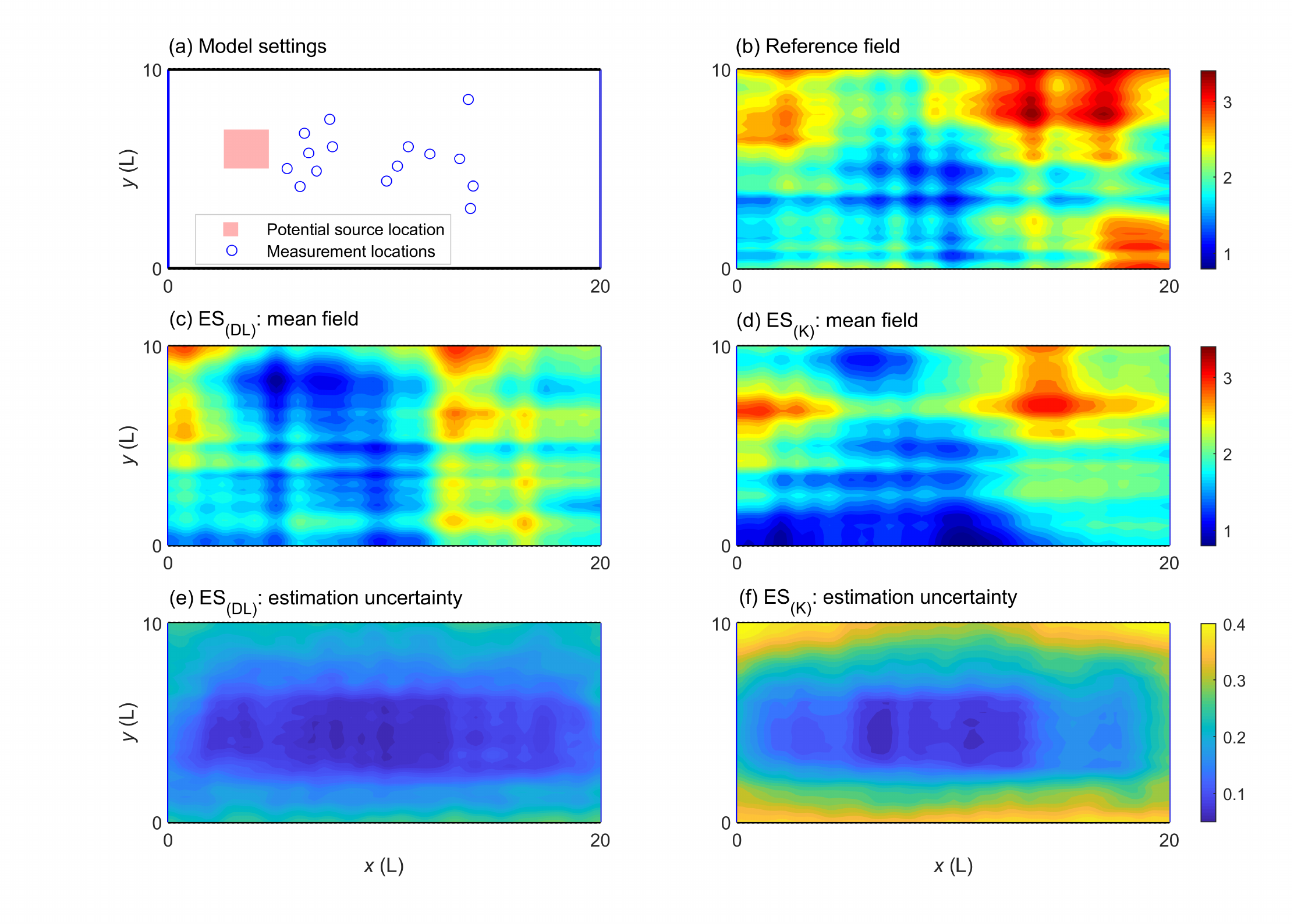}
	\caption{(a) Schematic overview of the flow domain in the first case study. Here, an unknown point source is located somewhere in the light red rectangle, and measurements of hydraulic head and solute concentration are collected at 15 monitoring wells denoted by the blue circles. (b) The reference log-conductivity ($\mathcal{Y}$) field. (c-f) Mean estimations of the $\mathcal{Y}$ field and the associated standard deviation fields obtained by ES$_\text{(DL)}$ and ES$_\text{(K)}$, respectively.}
	\label{fig: ex1_fields}
\end{figure}

In this section, we aim to demonstrate that when the variables involved are close to be Gaussian-distributed, ES$_\text{(DL)}$ can obtain very similar estimation of unknown model parameters to ES$_\text{(K)}$. Here, we consider an inverse problem where the parameters describing the hydraulic conductivity field and an unknown contaminant source are to be inferred from measurements of hydraulic head and solute concentration \cite{zhang2015,zhang2020}.

In this case, steady-state groundwater flow and transient solute transport are simulated in a two-dimensional (2D), heterogeneous, and confined aquifer. As shown in Figure \ref{fig: ex1_fields}a, the flow domain is 20 (L)$\times$10 (L) (in units of length), and discretized into 81$\times$41 grids in the numerical model. The left and right sides of the domain are prescribed by constant-head conditions of 12 (L) and 11 (L), respectively, while the upper and lower boundaries are impervious. At the initial time, the hydraulic head is 11 (L) everywhere in the domain, except for the left boundary. The hydraulic conductivity ($\mathcal{K}$) is heterogeneous and isotropic, and its logarithmic form, $\mathcal{Y}=\text{log}~\mathcal{K}$, is Gaussian-distributed and spatially correlated according to the following covariance function,
\begin{linenomath*}
	\begin{equation}
	C_{\mathcal{Y}}\left(x_{1}, y_{1} ; x_{2}, y_{2}\right)=\sigma_{\mathcal{Y}}^{2} \exp \left(-\frac{\left|x_{1}-x_{2}\right|}{\lambda_{x}}-\frac{\left|y_{1}-y_{2}\right|}{\lambda_{y}}\right),
	\label{eq: correlation}
	\end{equation}
\end{linenomath*}
where $\{x_1,y_1\}$ and $\{x_2,y_2\}$ are two arbitrary locations in the domain, $\sigma_{\mathcal{Y}}^2$ is the variance of the $\mathcal{Y}$ field, and $\lambda_{x}$ and $\lambda_{y}$ are the correlation lengths in the horizontal ($x$) and vertical ($y$) direction, respectively. The reference, or ``true" log-conductivity field is depicted in Figure \ref{fig: ex1_fields}b. With the above model settings, we can obtain steady-state hydraulic head, $h$ (L), by solving
\begin{linenomath*}
	\begin{equation}
	\frac{\partial}{\partial x_{i}}\left(\mathcal{K}_{i} \frac{\partial h}{\partial x_{i}}\right)=0,
	\label{eq: darcy}
	\end{equation}
\end{linenomath*}
and obtain the pore water velocity, $v_i$ (LT$^{-1}$), by solving 
\begin{linenomath*}
	\begin{equation}
	v_{i}=-\frac{\mathcal{K}_{i}}{\theta} \frac{\partial h}{\partial x_{i}},
	\label{eq: velocity}
	\end{equation}
\end{linenomath*}
numerically with MODFLOW \cite{harbaugh2000}. Here, $\theta$ (-) is the aquifer porosity, and the subscript $i$ denotes the coordinate axis ($i=1$ is for the $x$ direction, and $i=2$ is for the $y$ direction). 

In the flow domain, there is a point source that releases some non-reactive contaminant to the downstream. The contaminant source is located somewhere in the light red rectangular zone in Figure \ref{fig: ex1_fields}a. Its release strength varies with time and is characterized by a step function composed of six mass-loading rates, that is, $S_k$ (MT$^{-1}$) from $k$ (T) to $k+1$ (T), $k=1,...,6$. By numerically solving the following advection-dispersion equation,
\begin{linenomath*}
	\begin{equation}
	\frac{\partial(\theta C)}{\partial t}=\frac{\partial}{\partial x_{i}}\left(\theta D_{i j} \frac{\partial C}{\partial x_{j}}\right)-\frac{\partial}{\partial x_{i}}\left(\theta v_{i} C\right)+q_{\mathrm{a}} C_{\mathrm{s}},
	\label{eq: advection-dispersion}
	\end{equation}
\end{linenomath*}
with MT3DMS \cite{zheng1999}, we can obtain the simulated concentrations, $C$ (ML$^{-3}$), at different times and places. Here, $t$ (T) is the time, $q_\text{a}$ (T$^{-1}$) is the volumetric flow rate per unit volume of the aquifer, $C_\text{s}$ (ML$^{-3}$) denotes the concentration of the contaminant source, and $D_{ij}$ (L$^2$T$^{-1}$) signifies the hydrodynamic dispersion tensor that is composed of
\begin{linenomath*}
	\begin{equation}
	\begin{aligned}
	&D_{11}=\frac{1}{\|\mathbf{v}\|}\left(\alpha_{\mathrm{L}} v_{1}^{2}+\alpha_{\mathrm{T}} v_{2}^{2}\right),\\
	&D_{22}=\frac{1}{\|\mathbf{v}\|}\left(\alpha_{\mathrm{L}} v_{2}^{2}+\alpha_{\mathrm{T}} v_{1}^{2}\right),\\
	&D_{12}=D_{21}=\frac{1}{\|\mathbf{v}\|}\left(\alpha_{\mathrm{L}}-\alpha_{\mathrm{T}}\right) v_{1} v_{2},
	\end{aligned}
	\label{eq: dispersion}
	\end{equation} 
\end{linenomath*}
where $\alpha_\text{L}$ and $\alpha_\text{T}$ (L) represent the longitudinal and transverse diversity, respectively, and $\|\mathbf{v}\|=\sqrt{v_1^2+v_2^2}$ is the magnitude of the velocity vector, $\textbf{v}$.

\begin{table}
\caption{Prior ranges and true values of the eight contaminant source parameters in the first numerical experiment.}
\centering
\begin{tabular}{c c c c c c c c c}
\hline
Parameter   & $x_\text{s}$ & $y_\text{s}$ & $S_1$ & $S_2$ & $S_3$ & $S_4$ & $S_5$ & $S_6$ \\
\hline
Prior range & [3-5]        & [4-6]        & [0-8] & [0-8] & [0-8] & [0-8] & [0-8] & [0-8]  \\
True value  & 3.52         & 4.44         & 5.69  & 7.88  & 6.31  & 1.49  & 6.87  & 5.55   \\
\hline
\label{tab: 1}
\end{tabular}
\end{table}

In this case, the uncertainty comes from the heterogeneous $\mathcal{Y}$ field and the unknown contaminant source. To reduce the dimensionality of the $\mathcal{Y}$ field, the truncated Karhunen-Lo\`{e}ve (KL) expansion \cite{zhang2004} is used to represent the $\mathcal{Y}$ field, that is,
\begin{linenomath*}
	\begin{equation}
	\widetilde{\mathcal{Y}}(\mathbf{x})=\mu_{\mathcal{Y}}+\sum_{n=1}^{N_{\text{KL}}} \sqrt{\tau_{n}} s_{n}(\mathbf{x}) \xi_{n},
	\label{eq: KL}
	\end{equation}
\end{linenomath*}
where $\textbf{x}=\{x,y\}$ is the location, $\mu_{\mathcal{Y}}$ denotes the mean of the $\mathcal{Y}$ field, $\tau_n$ and $s_n(\textbf{x})$ signify the eigenvalues and eigenfunctions of the covariance defined in equation (\ref{eq: correlation}), and $\xi_n\sim\mathcal{N}(0,1^2)$ represent the KL expansion terms, $n=1,...,N_{\text{KL}}$. Here, $N_{\text{KL}}=100$ KL terms are kept, which can preserve about 95\% of the total field variance, that is, $\sum_{n=1}^{N_{\text{KL}}}\tau_n/\sum_{n=1}^{\infty}\tau_n\approx0.95$. The contaminant source is parameterized by eight variables, that is, its location, $\{x_\text{s},y_\text{s}\}$, and time-varying source strengths, $\{S_1,...,S_6\}$. Prior distributions of the eight source parameters are uniform and bounded by the ranges as listed in Table \ref{tab: 1}. Thus, there are 108 unknown parameters to be estimated in this case, that is, $\textbf{m}=\{\xi_1,...,\xi_{100},x_\text{s},y_\text{s},S_1,...,S_6\}$. Other model parameters are obtained from experiments or geological surveys as $\sigma_\mathcal{Y}^2=1$, $\lambda_x=10$ (L), $\lambda_y=5$ (L), $\mu_{\mathcal{Y}}=2$, $\theta=0.25$ (-), $\alpha_{\text{L}}=0.3$ (L), and $\alpha_{\text{T}}=0.03$ (L), respectively. 

\begin{figure}
	\noindent\includegraphics[width=\textwidth]{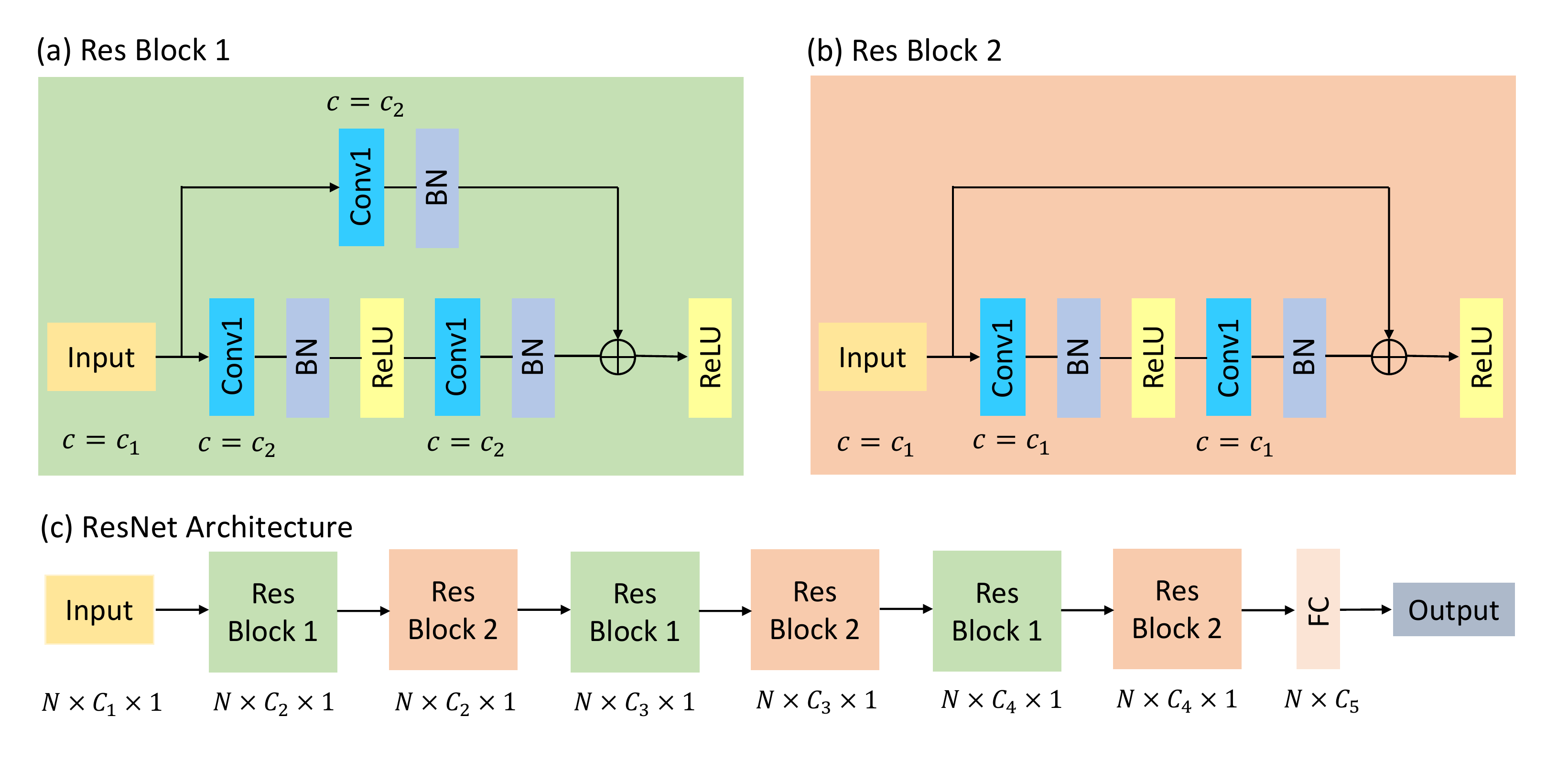}
	\caption{(a-b) Two residual (Res) blocks containing three kinds of layers, that is, 1D convolution (Conv1), batch normalization (BN), and rectified linear unit (ReLU) activation. \add[editor]{\protect Here, Conv1 applies convolutional kernels sliding along the 1D input; BN first normalizes its input to maintain zero mean and standard deviation, and then shifts and scales the data by a learnable offset and a learnable scale factor, respectively; ReLU sets all negative values in the input to zero;} In Res Block 1, the channels of input are scaled down before adding ($\bigoplus$) to the output, while in Res Block 2 the channels of input keep unchanged. (c) The overall architecture of the deep learning model (ResNet) used in the first case study. \add[editor]{\protect Here, FC means a fully-connected layer where the input data are first multiplied by a weighted matrix and then added to a bias vector.}}
	\label{fig: ex1_dl}
\end{figure}

To infer the 108 unknown model parameters, steady-state hydraulic heads, and transient solute concentrations at $t=\{4,5,...,12\}$ (T), are collected at 15 monitoring wells denoted by the blue circles in Figure \ref{fig: ex1_fields}a. The measurements are generated by running the integrated model (MODFLOW+MT3DMS) with the reference log-conductivity field (Figure \ref{fig: ex1_fields}b) and contaminant source parameters (the last row of Table \ref{tab: 1}), and adding independent normal random perturbations that satisfy $\upepsilon_h\sim\mathcal{N}(0,0.005^2)$ and $\upepsilon_C\sim\mathcal{N}(0,0.005^2)$ for hydraulic heads and solute concentrations, respectively. Then we implement the ES$_\text{(K)}$ and ES$_\text{(DL)}$ methods respectively to estimate the unknown model parameters conditioned on the measurement data. As the problem tested here is rather nonlinear, we perform multiple data assimilation ($N_\text{iter}=5$) in the two ES methods. At first, a same forecast ensemble ($N_\text{e}=500$) is generated from the prior parameter distribution for the two methods, $\textbf{M}^{(0)}=\{\textbf{m}_1^{(0)},...,\textbf{m}_{N_{\text{e}}}^{(0)}\}$, and the corresponding model outputs are calculated by running the numerical model, $\textbf{Y}^{(0)}=\{f(\textbf{m}_1^{(0)}),...,f(\textbf{m}_{N_\text{e}}^{(0)})\}$.

In each iteration of ES$_\text{(DL)}$, a same DL architecture as shown in Figure \ref{fig: ex1_dl} is adopted. The dimensions of inputs and outputs of the DL model are $150\times1$ and $108\times1$, respectively. In iteration $t~(1\le t\le N_\text{iter})$, to train the DL model, a set of data with $N=N_\text{e}(N_\text{e}-1)/2$ samples, that is, $\textbf{D}^{(t-1)}=\{\textbf{D}_\text{in}^{(t-1)}~\textbf{D}_\text{out}^{(t-1)}\}$, are generated from $\{\textbf{M}^{(t-1)}~\textbf{Y}^{(t-1)}\}$. Here, $\textbf{D}_\text{in}^{(t-1)}=\{f(\textbf{m}_i^{(t-1)})-f(\textbf{m}_j^{(t-1)})+\alpha_t\upepsilon_{ij}|i=1,...,N_\text{e}-1,~i<j\le N_\text{e}\}$ are the input data, $\textbf{D}_\text{out}^{(t-1)}=\{\textbf{m}_i^{(t-1)}-\textbf{m}_j^{(t-1)}|i=1,...,N_\text{e}-1,~i<j\le N_\text{e}\}$ are the output data, $\alpha_t$ is an inflation factor that can be conveniently set as $\sqrt{N_\text{iter}}$, and $\upepsilon_{ij}$ are random realizations of the measurement error. Here, to sufficiently extract features embedded in the training data, $\textbf{D}^{(t-1)}$, we employ the residual network (ResNet) proposed by \citeA{he2016}. To adapt to our data format, we replace the 2D convolution (suitable for image-like data) used in the original ResNet with one-dimensional (1D) convolution (i.e., Conv1, suitable for sequence-like data). It is noted here that the kernel size and stride in Conv1 are both set as 1, thus Conv1 works similarly to a fully-connected (FC) layer. In ResNet, we can build a very deep network without worrying about the trouble caused by gradient vanishing. As shown in Figure \ref{fig: ex1_dl}c, the overall architecture is composed of two kinds of residual blocks (Res Block 1 in Figure \ref{fig: ex1_dl}a and Res Block 2 in Figure \ref{fig: ex1_dl}b) and a FC layer. The inputs to these blocks are all vectors. The difference between Res Block 1 and Res Block 2 lies in that the number of channels in the former block is scaled down, while in the latter block the number of channels is unchanged. In this case, the numbers of channels are designed as $C_1=150$, $C_2=140$, $C_3=130$, $C_4=120$, and $C_5=108$, respectively. The Adam optimizer \cite{kingma2014} with a learning rate of $3\times 10^{-3}$ is utilized to train the network.

\begin{figure}
	\noindent\includegraphics[width=\textwidth]{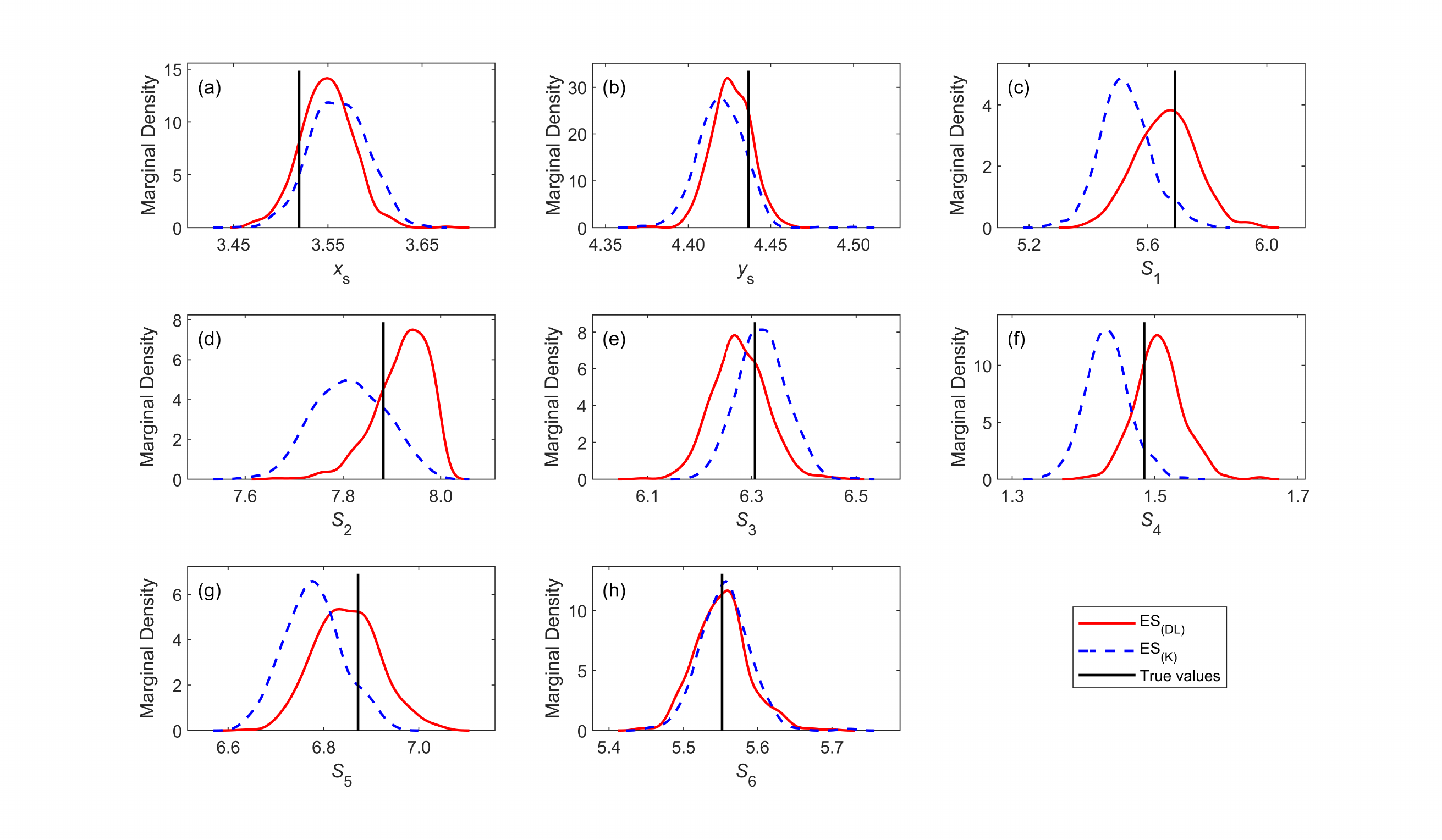}
	\caption{Marginal distributions of the eight contaminant source parameters estimated by ES$_\text{(DL)}$ (red lines) and ES$_\text{(K)}$ (blue dashed lines), respectively. The true values of the eight parameters are denoted by the vertical black lines.}
	\label{fig: ex1_ppdf}
\end{figure}

After five iterations, both ES$_\text{(K)}$ and ES$_\text{(DL)}$ can significantly improve our knowledge of the subsurface medium and contaminant source. As shown in Figures \ref{fig: ex1_fields}(c-d), both methods can reliably identify the regions with high and low values of (log) conductivity. Yet the two estimated mean fields tend to slightly underestimate the true $\mathcal{Y}$ values. The root-mean-square errors (RMSEs) between the estimated mean fields and the reference field (Figure \ref{fig: ex1_fields}b) are 0.4690 and 0.5147 for ES$_\text{(K)}$ and ES$_\text{(DL)}$, respectively. Figures \ref{fig: ex1_fields}(e-f) present the standard deviation (SD) fields associated with the mean estimates. It can be found that the area where monitoring wells have been installed exhibits smaller SD values, and the results from ES$_\text{(DL)}$ have smaller variations than that of ES$_\text{(K)}$. In Figure \ref{fig: ex1_ppdf}, we draw the marginal densities of the eight contaminant source parameters estimated from the updated ensemble in the last iteration, that is, $\textbf{M}^{(N_\text{iter})}=\{\textbf{m}_1^{(N_\text{iter})},...,\textbf{m}_{N_\text{e}}^{(N_\text{iter})}\}$. Here, we use red lines and blue dashed lines to represent the results from ES$_\text{(DL)}$ and ES$_\text{(K)}$, respectively. Compared to the prior ranges as listed in Table \ref{tab: 1}, the ranges covered by the marginal densities are much narrower, which indicates a substantial reduction of uncertainty in our belief about the model parameters. Moreover, the true parameter values (vertical black lines) generally locate near the centers of the marginal density curves, which indicates the accuracy of the estimation results. For the eight contaminant source parameters, we calculate the root-mean-square relative errors (RMSREs) between the mean estimates and the true parameter values. The corresponding RMSRE values for ES$_\text{(K)}$ and ES$_\text{(DL)}$ are 0.0177 and 0.0066, respectively. 

From the above results, it is found that ES$_\text{(K)}$ performs slightly better at characterizing the log-conductivity field, while ES$_\text{(DL)}$ can more accurately identify the contaminant source parameters. Overall, the two ES methods can obtain reliable and comparable estimations of the log-conductivity field and unknown contaminant source parameters. If a more diverse and larger measurement dataset is collected, and/or a more suitable DL architecture is designed, the ES$_\text{(DL)}$ method should be able to produce better results.

\subsection{Example 2: A Non-Gaussian Case }

In the previous section, we have tested a case where the distributions of concerned variables are near multi-Gaussian, and ES$_\text{(DL)}$ can produce similar results as ES$_\text{(K)}$. Nevertheless, in subsurface characterization, much research has shown that when the parameter field of interest does not follow a multi-Gaussian distribution, the direct use of a Kalman-based DA method, for example, EnKF or ES$_\text{(K)}$, cannot produce satisfactory results \cite{cao2018,chang2010,xu2016,zhou2011}. Below we will test such a case where sparse measurements of hydraulic head are used to characterize a non-Gaussian conductivity field, and the performances of ES$_\text{(DL)}$ and ES$_\text{(K)}$ are compared.

\begin{figure}
	\noindent\includegraphics[width=\textwidth]{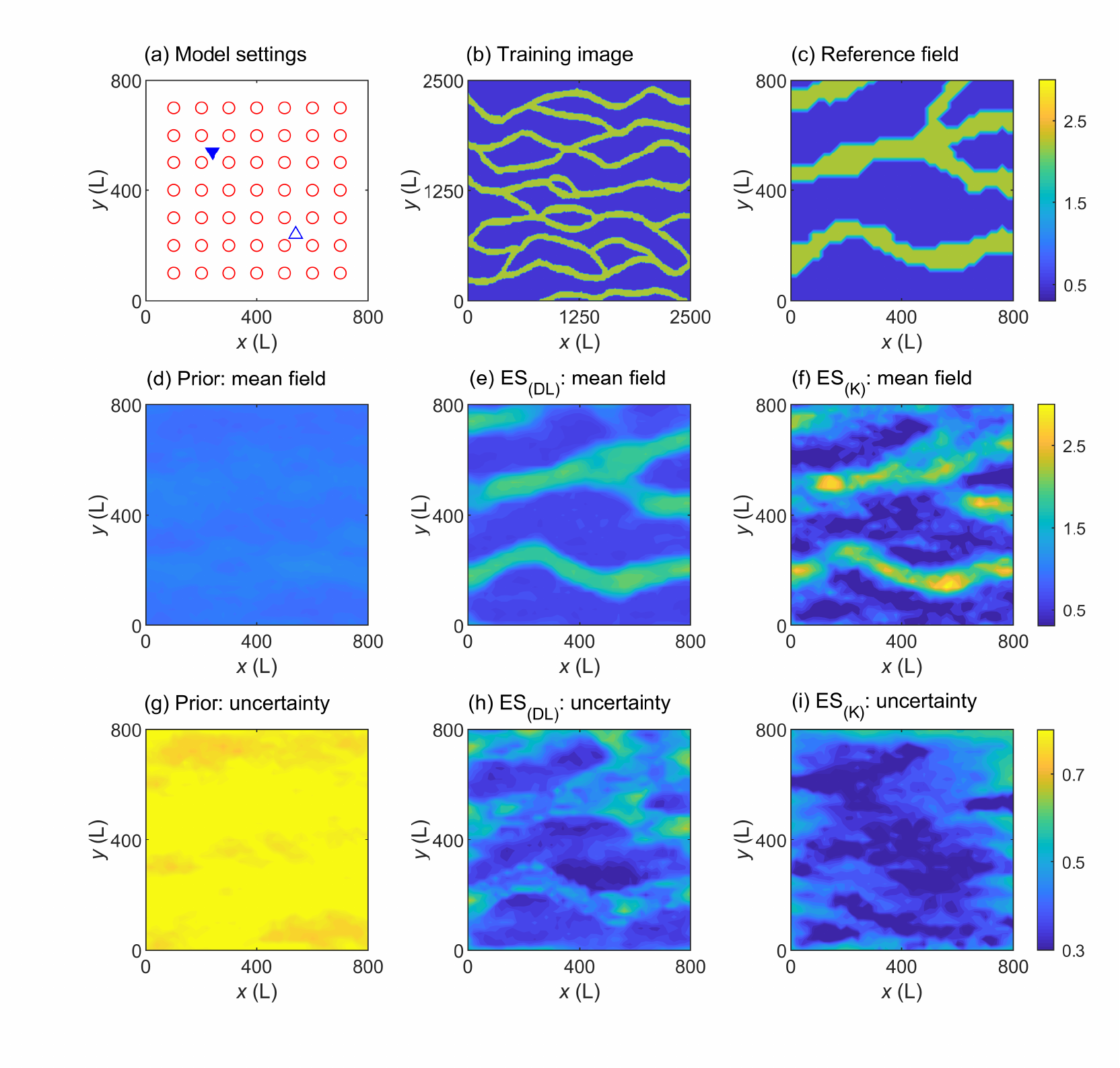}
	\caption{(a) Schematic overview of the flow domain in the second case study. Here, the injection well is denoted by the blue down-pointing triangle, the pumping well is represented by the blue up-pointing triangle, and measurements of hydraulic head are collected at wells signified by the 7$\times$7 red circles. (b) The training image used in the direct sampling method. (c) The reference conductivity ($\mathcal{K}$) field. (d-h) Mean estimations of the $\mathcal{K}$ field and their associated standard deviations from the prior ensemble and the updated ensembles obtained by ES$_\text{(DL)}$ and ES$_\text{(K)}$, respectively.}
	\label{fig :ex2_fields}
\end{figure}

Here, we consider transient water flow in a 2D, confined, and channelized aquifer. The size of the domain is $l_x=l_y=800$ (L) in units of length in the $x$ and $y$ direction. This square domain is uniformly discretized into 41$\times$41 grids. In the flow field, impervious condition is prescribed at both the upper and lower boundaries, and constant heads of 202 (L) and 198 (L) are imposed at the left and right sides, respectively. At the initial time, the hydraulic head is 198 (L) across the domain except for the left boundary. To enhance water flow in the subsurface medium, an injection well (the blue down-pointing triangle in Figure \ref{fig :ex2_fields}a) with a rate of 150 ($\text{L}^3\text{T}^{-1}$) and a pumping well (the blue up-pointing triangle in Figure \ref{fig :ex2_fields}a) with a rate of -150 ($\text{L}^3\text{T}^{-1}$) are installed. In the channelized field, there are two kinds of materials: one with a low conductivity value of $\mathcal{K}_1=0.5~(\text{LT}^{-1})$, and another with a higher value of $\mathcal{K}_2=2.3~(\text{LT}^{-1}$). The reference $\mathcal{K}$ field (Figure \ref{fig :ex2_fields}c) is generated from a training image (Figure \ref{fig :ex2_fields}b) using the direct sampling (DS) method proposed by \citeA{mariethoz2010}. Here, when applying the DS method, no direct observation of $\mathcal{K}$ is used for conditioning. The DS method is computationally efficient, and it has the ability to handle both continuous and categorical variables with complex patterns. Thus, it is adopted in this case to perform multiple-point statistics simulations to generate the reference, as well as random realizations of non-Gaussian $\mathcal{K}$ field. Details of the DS method can be found in \cite{mariethoz2010,meerschman2013}. With the above model settings, one can obtain transient hydraulic heads at different locations, $h(\mathbf{x},t)~(\text{L})$, by solving
\begin{linenomath*}
	\begin{equation}
	S_{\mathrm{s}} \frac{\partial h(\mathbf{x}, t)}{\partial t}+\nabla \cdot \mathbf{q}(\mathbf{x}, t)=g(\mathbf{x}, t)
	\label{eq: transient}
	\end{equation}
\end{linenomath*}
numerically with MODFLOW \cite{harbaugh2000}. In equation (\ref{eq: transient}), $S_{\mathrm{s}}~(\text{L}^{-1})$ is the specific storage, $\mathbf{x}=\{x,y\}~(\text{L})$ is the location, $t~(\text{T})$ is the time, $\mathbf{q}(\mathbf{x}, t)=-\mathcal{K}(\mathbf{x})\nabla h(\mathbf{x},t)$ is the flux, $\mathcal{K}(\mathbf{x})~(\text{LT}^{-1})$ is the conductivity value at location $\mathbf{x}$, and $g(\mathbf{x}, t)~(\text{T}^{-1})$ is the source (or sink) term of water. Here, the total simulation time is 18 (T), and $S_{\mathrm{s}}$ is a deterministic constant of 0.0001 (L$^{-1}$).

To infer the $\mathcal{K}$ field, we collect measurements of hydraulic head at 49 wells denoted by the red circles in Figure \ref{fig :ex2_fields}a, every 0.6 (T) from $t=0$ (T) to $t=6$ (T). The measurements are generated by running the numerical model with the reference $\mathcal{K}$ field (Figure \ref{fig :ex2_fields}c) and adding perturbations that fit $ \upepsilon\sim\mathcal{N}_n(0,0.01^2)$. For both the ES$_\text{(DL)}$ and ES$_\text{(K)}$ methods, a same set of $N_\text{e}=499$ prior random realizations of channelized field, that is, $\textbf{M}^{(0)}=\{\textbf{m}_1^{(0)},...,\textbf{m}_{N_\text{e}}^{(0)}\}$, are generated using the DS method based on the training image (Figure \ref{fig :ex2_fields}b). By averaging these $N_\text{e}$ realizations, we can obtain a rather uniform prior mean field (Figure \ref{fig :ex2_fields}d) with grid values close to 0.98, the mean value of the training image (averaged over each grid). The associated standard deviation field (Figure \ref{fig :ex2_fields}g) also exhibits a small spatial variability, and has values close the standard deviation (0.80) of the training image. Through running the numerical model, we can obtain the corresponding model outputs, that is, $\textbf{Y}^{(0)}=\{f(\textbf{m}_1^{(0)}),...,f(\textbf{m}_{N_\text{e}}^{(0)})\}$.

\begin{figure}
	\noindent\includegraphics[width=\textwidth]{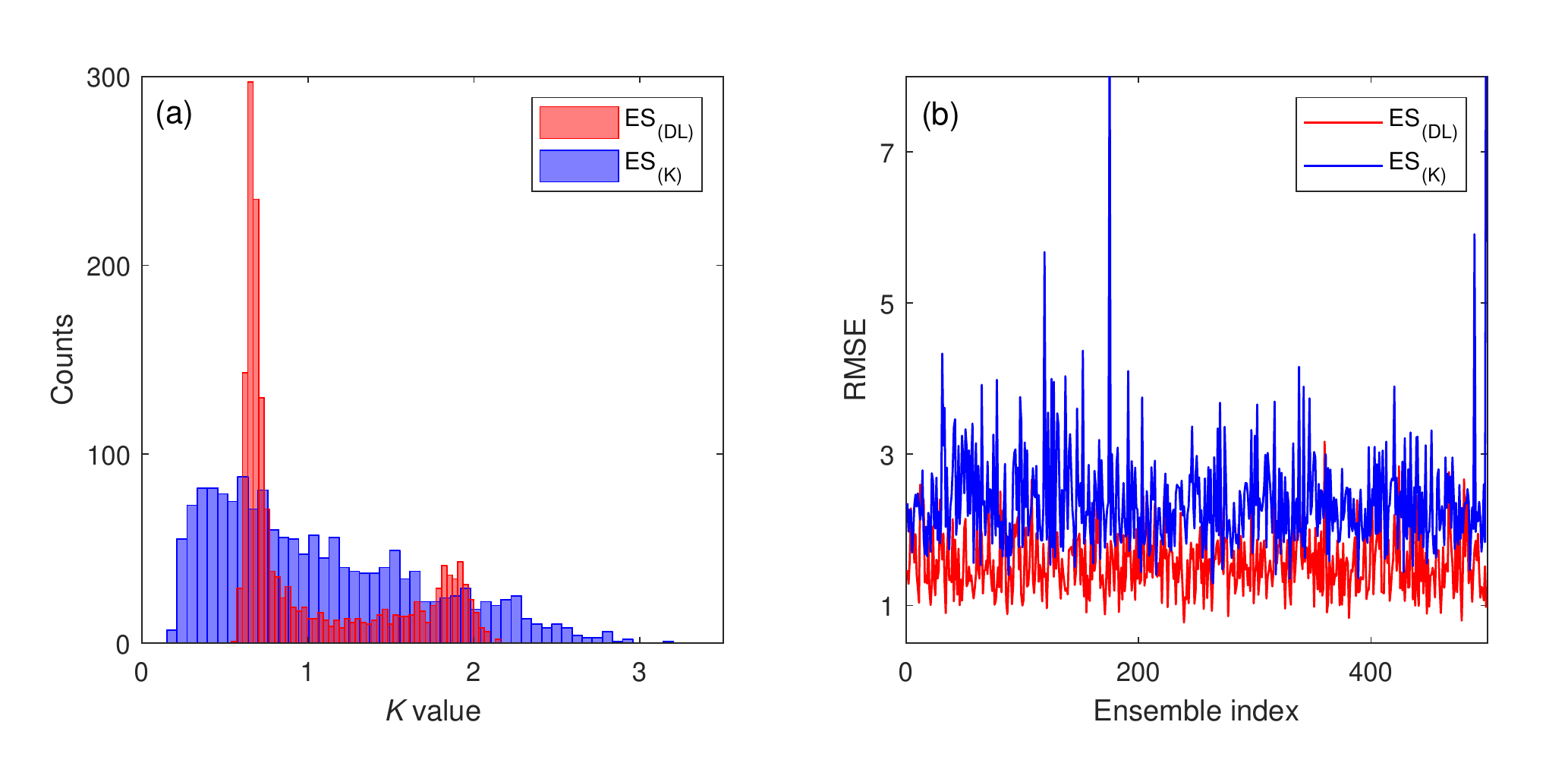}
	\caption{(a) Histograms of the mean conductivity ($\mathcal{K}$) fields estimated by ES$_\text{(DL)}$ and ES$_\text{(K)}$, respectively. (b) Root-mean-square error (RMSE) of the MODFLOW simulated and observed hydraulic heads for each sample in the updated ensembles of ES$_\text{(DL)}$ and ES$_\text{(K)}$. Color coding in red and blue differentiates between the results from ES$_\text{(DL)}$ and ES$_\text{(K)}$.}
	\label{fig :ex2_hist_rmse}
\end{figure}

Figure \ref{fig :ex2_fields}f presents the mean $\mathcal{K}$ field estimated by the ES$_\text{(K)}$ method. As the problem considered here is rather linear, it is not necessary to perform multiple data assimilation, that is, here we set $N_\text{iter}=1$. It is obvious that ES$_\text{(K)}$ can capture some patterns of the true $\mathcal{K}$ field through assimilating indirect measurements of transient hydraulic head. Moreover, the SD field calculated from the updated ensemble (Figure \ref{fig :ex2_fields}i) has much smaller values than the prior SD field (Figure \ref{fig :ex2_fields}g). Nevertheless, the connectivity pattern of the mean $\mathcal{K}$ field is underestimated. In Figure \ref{fig :ex2_hist_rmse}a, we draw the histogram of the mean $\mathcal{K}$ field of ES$_\text{(K)}$ (blue bars). As the true $\mathcal{K}$ field only has two distinct materials, ideally, the histogram should be bimodal. However, ES$_\text{(K)}$ fails to recover this bimodality.

\begin{figure}
	\noindent\includegraphics[width=\textwidth]{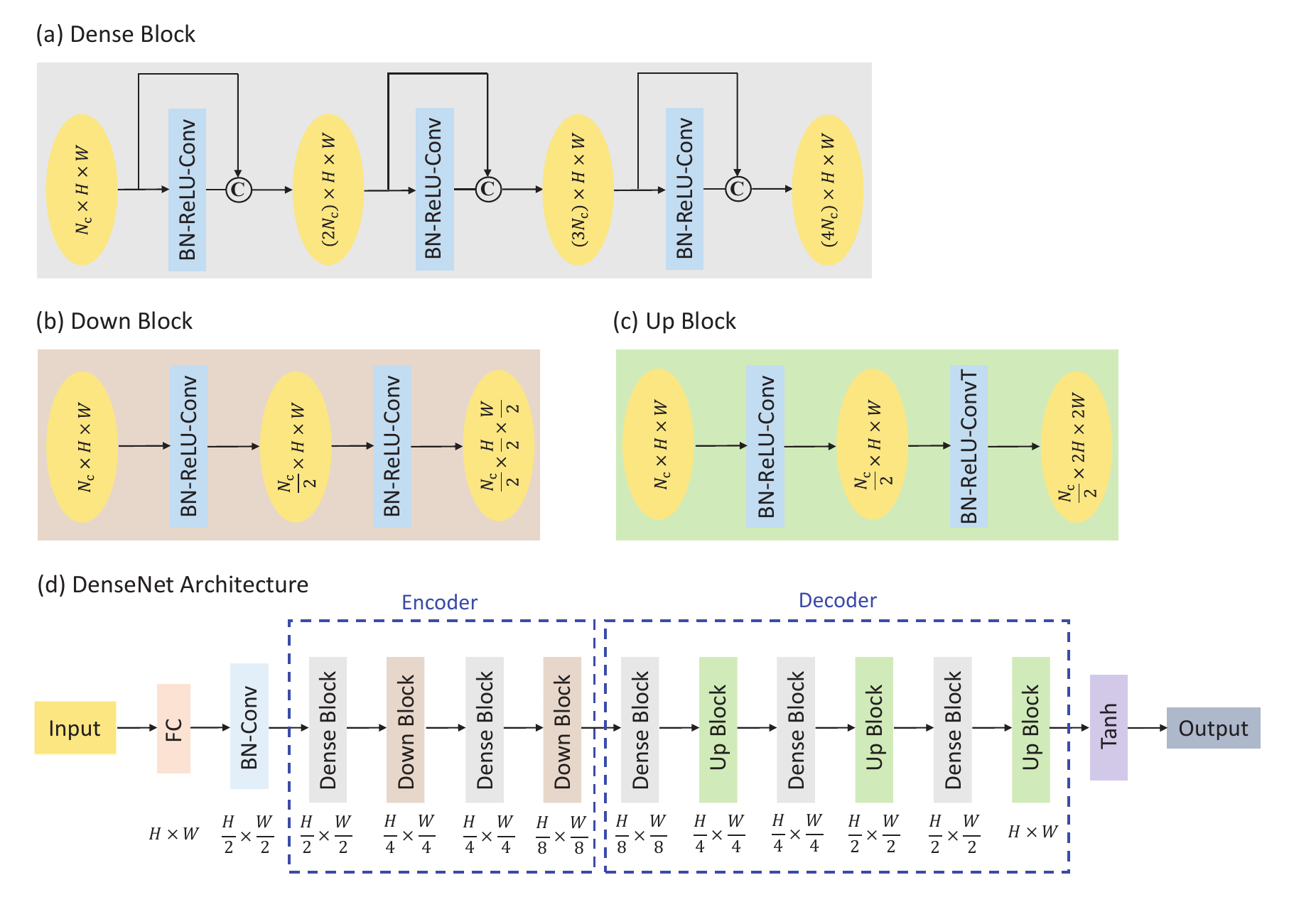}
	\caption{(a) A dense block with three layers, each of which performs three operations sequentially, that is, batch normalization (BN), rectified linear unit (ReLU) activation, and 2D convolution (Conv). Here, the feature maps (i.e., the outputs) of each layer will be concatenated ($\textcopyright$) with their corresponding inputs. (b) A transition block used to reduce both the number of channels and the size of feature maps. (c) A transition block used to reduce the number of channels but increase the size of feature maps by using transpose convolution (ConvT). (d) The overall architecture of the deep learning model (DenseNet) used in the second case study. \add[editor]{\protect Here, FC means a fully-connected layer where the input data are first multiplied by a weighted matrix and then added to a bias vector; BN first normalizes its input to maintain zero mean and standard deviation, and then shifts and scales the data by a learnable offset and a learnable scale factor, respectively; Conv applies convolutional kernels sliding along the 2D input to extract feature maps; ConvT is a transposed Conv layer that upsamples the feature maps (i.e., make them larger); ReLU sets all negative values in the input to zero; Tanh applies the tanh function on the input.}}
	\label{fig: ex2_dl}
\end{figure} 

Then we apply the ES$_\text{(DL)}$ method to estimate the $\mathcal{K}$ field. 
In this case, the input and output dimensions of the DL model are $490\times1$ (i.e., the number of measurement data) and $41\times41$ (i.e., the number of model grids), respectively. Without extra evaluations of the numerical model, a set of training data, $\textbf{D}^{(0)}=\{\textbf{D}_\text{in}^{(0)}~\textbf{D}_\text{out}^{(0)}\}$, can be generated from the forecast ensemble, $\{\textbf{M}^{(0)}~\textbf{Y}^{(0)}\}$. Here, $\textbf{D}_\text{in}^{(0)}=\{f(\textbf{m}_i^{(0)})-f(\textbf{m}_j^{(0)})+\upepsilon_{ij}|i=1,...,N_\text{e}-1,~i<j\le N_\text{e}\}$ are the input data, $\textbf{D}_\text{out}^{(0)}=\{\textbf{m}_i^{(0)}-\textbf{m}_j^{(0)}|i=1,...,N_\text{e}-1,~i<j\le N_\text{e}\}$ are the output data, and $\upepsilon_{ij}$ are random realizations of the measurement error. A better analysis scheme is expected to be learned for ES from the training data, $\textbf{D}^{(0)}$, with an adequate DL model. Considering the fact that the monitoring wells (the 7$\times$7 red circles in Figure \ref{fig :ex2_fields}a) are uniformly distributed in the flow domain, here we suppose they have some spatial connectivity. Besides, the outputs corresponding to the update to the conductivity field can be naturally seen as an image. Therefore, the mapping from the innovation $\Delta\textbf{y}$ to the update $\Delta\textbf{m}$ can be transformed to an image-to-image task, and 2D convolution is adopted to process the spatial features. When designing the DL model, we consider the popular DenseNet architecture proposed by \citeA{huang2017}. The DenseNet architecture enables each layer in the network to connect with any previous layer and realizes feature re-utilization, thus it can reduce the redundancy of the training parameters and improve the efficiency remarkably. As shown in Figure \ref{fig: ex2_dl}d, the overall architecture of DenseNet is composed of an encoder, a decoder, and some other necessary layers. The encoder aims to discover low-dimensional embeddings of the input image, while the decoder maps the embeddings to the output image. The encoder and decoder are composed of three kinds of basis blocks as shown in Figures \ref{fig: ex2_dl}(a-c). The dense block (Figure \ref{fig: ex2_dl}a) is a concatenation of previous feature maps. In this case, the dense block contains three layers that perform batch normalization (BN), rectified linear unit (ReLU) activation, and 2D convolution (Conv) sequentially. Specifically, the Conv operation is realized by a kernel of size $3 \times 3$ with a stride of 1 and a padding of 1. As a result, the dense block produces outputs with channels four times as many as that in the inputs, while the size of feature maps keeps unchanged. To adjust the size of feature maps, two kinds of transition blocks are further employed, that is, the down block (Figure \ref{fig: ex2_dl}b) and the up block (Figure \ref{fig: ex2_dl}c). The two blocks both contain two convolution layers. The first convolution layers in the two blocks work in the same way that uses a $1 \times 1$ kernel with a stride of 1. However, their second convolution layers are implemented differently in that the down block utilizes a $3\times 3$ kernel with a stride of 2 and a padding of 1 for downsampling, while the up block uses a kernel of the same size to perform transpose convolution (ConvT) for upsampling. Overall, as shown in Figure \ref{fig: ex2_dl}d, the input data are first processed by a fully-connected (FC) layer to obtain image-like data of size $H \times W$, and then go to the BN-Conv layers, the encoder-decoder blocks, and finally a Tanh activation layer to produce the output data. In this case, $H=W=41$, the number of channels after the FC layer is $N_\text{c}=40$, and after the decoder block the number of channels is $N_\text{c}=1$. The Adam optimizer with a learning rate of $1\times10^{-3}$ is utilized to train the network.

As shown in Figure \ref{fig :ex2_fields}e, the mean conductivity field estimated by ES$_\text{(DL)}$ better resembles the reference $\mathcal{K}$ field. Besides, the RMSE value between the mean estimate from ES$_\text{(DL)}$ and the reference $\mathcal{K}$ field is 0.5010, which is smaller than the RMSE value of 0.5607 from ES$_\text{(K)}$. Although the SD field of ES$_\text{(DL)}$ has slightly larger values than the SD field of ES$_\text{(K)}$, the channelized features are better revealed in Figure \ref{fig :ex2_fields}h. Moreover, the histogram of the mean $\mathcal{K}$ field of ES$_\text{(DL)}$ (red bars, Figure \ref{fig :ex2_hist_rmse}a) can clearly recover the bimodality of the channelized field, although the update of $\mathcal{K}_1$ is slightly overestimated, while the update of $\mathcal{K}_2$ is slightly underestimated. Thus, we believe that ES$_\text{(DL)}$ can better handle non-Gaussian parameter field than ES$_\text{(K)}$. In Figure \ref{fig :ex2_hist_rmse}b, we draw the RMSE of the MODFLOW simulated and observed hydraulic heads for each sample in the updated ensembles of ES$_\text{(DL)}$ (red line) and ES$_\text{(K)}$ (blue line). It again demonstrates the superiority of ES$_\text{(DL)}$ to ES$_\text{(K)}$. If a better DL architecture is designed, a more accurate update of the parameter field can be obtained by the ES$_\text{(DL)}$ method.

\section{Discussions and Conclusions}
Due to their efficiency and robustness, EnKF and its variants have been used in various research fields of geosciences to reduce the uncertainty of the system-of-interest. When one's purpose is parameter estimation, for example, in subsurface characterization, ES can be adopted as a feasible method. Nevertheless, when the distributions of involved variables are non-Gaussian, performances of these Kalman-based DA methods will deteriorate. To enable proper applications of these methods, existing strategies mainly transform non-Gaussian variables to be normally distributed \cite{zhou2011,chang2010,canchumuni2019}, or use another method, for example, clustering analysis, or a more general DA method like particle filter, to handle non-Gaussianity \cite{cao2018,sun2009,mandel2009}.

Alternatively, we propose in this work to use DL to reformulate the analysis scheme of ES to gain an improved performance. In this new method, that is, ES$_\text{(DL)}$, we first generate a high volume of training data from a relatively small-sized forecast ensemble. Possible non-Gaussian features in model parameters and observations are incorporated in the training data and captured by an adequate DL model. Then we use this DL-based formulation to update the forecast ensemble to reduce the uncertainty of model parameters. For highly nonlinear problems, an iterative application of ES is needed, for example, using the multiple data assimilation scheme formulated by \citeA{emerick2013}. To demonstrate the performance of the proposed method, two cases of subsurface characterization are tested against the traditional ES method using the Kalman formula, that is, ES$_\text{(K)}$. In the first case study, using measurements of hydraulic head and solute concentration, we aim to simultaneously identify the location and release history of a point contaminant source, as well as the heterogeneous log-conductivity field. Here, there are 108 unknown parameters to be estimated, whose distributions are all close to Gaussian. With the same number of numerical model evaluations, the ES$_\text{(DL)}$ method produces comparable results to those from ES$_\text{(K)}$. In the second case study, a channelized conductivity field parameterized by 1681 variables is to be estimated from observations of transient hydraulic head. Simulation results clearly indicate that, in this non-Gaussian case, ES$_\text{(DL)}$ is superior to ES$_\text{(K)}$. 

The general applicability of ES$_\text{(DL)}$ comes from the powerful ability of DL in extracting complex (including non-Gaussian) features and learning nonlinear relationships automatically from data. The DL architecture is very flexible and can be adapted to a wide range of problems. Without running a large number of system models, one can create massive amounts of training data and feed them to the DL model. Another merit of DL is that it can perform massively parallel computations on GPUs. Thus, the ES$_\text{(DL)}$ method can be possibly applied to large-scale DA problems. Nevertheless, limitations of the proposed method do exist. First, the choice of a DL architecture for ES$_\text{(DL)}$ is relatively subjective (although flexible), and its outputs are difficult to comprehend. It is true that one can learn from literature to configure an adequate DL model, and different DL models can possibly all produce satisfactory results. However, there is no standard guideline to determine the optimal DL architecture for a specific problem. On the contrary, the Kalman formula used in ES$_\text{(K)}$ can be expressed explicitly, and it is optimal at least for linear, Gaussian cases. From the theoretical perspective, ES$_\text{(K)}$ is more elegant than ES$_\text{(DL)}$. Second, although the ES$_\text{(DL)}$ method requires a same number of system model evaluations as ES$_\text{(K)}$, the training of a DL model can be time consuming, especially when GPU devices are not available. Moreover, in this work, we only apply the new idea in parameter estimation problems. We believe that one can easily extend the DL-based idea to state estimations for real-time forecasting. In this case, it is natural to consider using recurrent neural networks (e.g., the famous long short-term memory network) to implement the DL-based idea in sequential DA problems. Recently, model structural uncertainty has been accounted for in the application of various iterative ES methods \cite{evensen2019}, which is important to prevent unphysical updates. When stochastic model errors are considered, one can rewrite equation (\ref{eq: numerical_model}) in the following way,
\begin{linenomath*}
	\begin{equation}
	\widetilde{\textbf{y}}=f(\textbf{m},
	\textbf{q})+\upepsilon,
	\label{eq: err_model}
	\end{equation}
\end{linenomath*}
where $\textbf{q}$ represent the model errors, and one simple form of $f(\textbf{m},
\textbf{q})$ can be chosen as $f(\textbf{m},
\textbf{q})=f(\textbf{m})+\textbf{q}$. In ensemble smoother, \citeA{evensen2019} proposed to update each prior sample of $\textbf{m}$ and $\textbf{q}$ as follows,
\begin{linenomath*}
	\begin{equation}
	\begin{split}
	\begin{array}{c}
	\mathbf{m}_{i}^{\mathrm{(1)}}=\mathbf{m}_{i}^{\mathrm{(0)}}+\overline{\mathbf{C}}_{\mathrm{MY}}^{\mathrm{(0)}}\left(\widetilde{\mathbf{C}}_{\mathrm{YY}}^{\mathrm{(0)}}+\mathbf{R}\right)^{-1}\left[\widetilde{\mathbf{y}}+\upepsilon_{i}-f\left(\mathbf{m}_{i}^{\mathrm{(0)}},\mathbf{q}_{i}^{\mathrm{(0)}}\right)\right], \\
	\mathbf{q}_{i}^{\mathrm{(1)}}=\mathbf{q}_{i}^{\mathrm{(0)}}+\overline{\mathbf{C}}_{\mathrm{QY}}^{\mathrm{(0)}}\left(\widetilde{\mathbf{C}}_{\mathrm{YY}}^{\mathrm{(0)}}+\mathbf{R}\right)^{-1}\left[\widetilde{\mathbf{y}}+\upepsilon_{i}-f\left(\mathbf{m}_{i}^{\mathrm{(0)}},\mathbf{q}_{i}^{\mathrm{(0)}}\right)\right], \\
	\end{array}
	\end{split}
	\label{eq: ES_error}
	\end{equation}
\end{linenomath*}
where $\overline{\mathbf{C}}_{\mathrm{MY}}^{\mathrm{(0)}}$, $\widetilde{\mathbf{C}}_{\mathrm{YY}}^{\mathrm{(0)}}$, and $\overline{\mathbf{C}}_{\mathrm{QY}}^{\mathrm{(0)}}$ are sample covariances calculated from the prior ensembles of model parameters, simulation outputs, and errors. Similarly, one can use DL to derive two new mappings to replace the two linear mappings defined by $\overline{\mathbf{C}}_{\mathrm{MY}}^{\mathrm{(0)}}\left(\widetilde{\mathbf{C}}_{\mathrm{YY}}^{\mathrm{(0)}}+\mathbf{R}\right)^{-1}$ and $\overline{\mathbf{C}}_{\mathrm{QY}}^{\mathrm{(0)}}\left(\widetilde{\mathbf{C}}_{\mathrm{YY}}^{\mathrm{(0)}}+\mathbf{R}\right)^{-1}$. For nonlinear problems, some iterative form of ES$_\text{(DL)}$ can be implemented. In future works, these ideas will be tested.

\acknowledgments
Computer codes and data used are available at \url{https://www.researchgate.net/publication/339447370_Using_Deep_Learning_to_Improve_Ensemble_Smoother}.\\
This work is supported by the National Key Research
and Development Program of China (grant 2018YFC1800503), and National Natural Science Foundation of China (grants 41807006 and 41771254). The authors would also like to thank Gregoire Mariethoz from University of Lausanne, Switzerland for providing the MATLAB codes of the direct sampling method.

\bibliography{myref}
\clearpage

\end{document}